\documentclass[12pt]{article}
\usepackage{amssymb}
\usepackage{epsfig}

\newcommand{\mb}[1]{ \mbox{\boldmath$#1$} }
\newcommand{\norm}[1]{\|#1\|}

\newcommand{\beq}{\begin{equation}}
\newcommand{\eeq}{\end{equation}}
\newcommand{\beqr}{\begin{eqnarray}}
\newcommand{\eeqr}{\end{eqnarray}}
\newcommand{\p}{\partial}

\newcommand{\ve}{\varepsilon}
\newcommand{\x}{\mb{x}}
\newcommand{\y}{\mb{y}}

\begin{document}

\title{Diffusion Maps, Spectral Clustering and Reaction Coordinates of Dynamical Systems}
\author{Boaz Nadler\footnotemark[2], St\'{e}phane Lafon\footnotemark[2], Ronald R. Coifman\footnotemark[2],
Ioannis G. Kevrekidis\footnotemark[3]
}

\maketitle
\renewcommand{\thefootnote}{\fnsymbol{footnote}}
\footnotetext[2]{Department of Mathematics, Yale University, New Haven, CT 06520.}
\footnotetext[3]{Chemical Engineering and PACM, Princeton University, Princeton, NJ 08544. }
\renewcommand{\thefootnote}{\arabic{footnote}}

\begin{abstract}

A central problem in data analysis is the low dimensional
representation of high dimensional data, and the concise
description of its underlying geometry and density. In the
analysis of large scale simulations of complex dynamical systems,
where the notion of time evolution comes into play, important
problems are the identification of slow variables and dynamically
meaningful reaction coordinates that capture the long time
evolution of the system. In this paper we provide a unifying view
of these apparently different tasks, by considering a family of
{\em diffusion maps}, defined as the embedding of complex (high
dimensional) data onto a low dimensional Euclidian space, via the
eigenvectors of suitably defined random walks defined on the given
datasets. Assuming that the data is randomly sampled from an
underlying general probability distribution $p(\x)=e^{-U(\x)}$, we
show that as the number of samples goes to infinity, the
eigenvectors of each diffusion map converge to the eigenfunctions
of a corresponding differential operator defined on the support of
the probability distribution. Different normalizations of the
Markov chain on the graph lead to different limiting differential
operators. For example, the normalized graph Laplacian leads to a
backward Fokker-Planck operator with an underlying potential of $2
U(\x)$, best suited for spectral clustering. A specific
anisotropic normalization of the random walk leads to the backward
Fokker-Planck operator with the potential $U(\x)$, best suited for
the analysis of the long time asymptotics of high dimensional
stochastic systems governed by a stochastic differential equation
with the same potential $U(\x)$. Finally, yet another
normalization leads to the eigenfunctions of the Laplace-Beltrami
(heat) operator on the manifold in which the data resides, best
suited for the analysis of the geometry of the dataset, regardless
of its possibly non-uniform density.
\end{abstract}

\section{Introduction}\label{s:intro}

Analysis of complex high dimensional data is an exploding area of
research, with applications in diverse fields, such as machine
learning, statistical data analysis, bio-informatics, meteorology,
chemistry and physics. In the first three application fields, the
underlying assumption is that the data is sampled from some
unknown probability distribution, typically without any notion of
time or correlation between consecutive samples. Important tasks
are dimensionality reduction, e.g., the representation of the high
dimensional data with only a few coordinates, and the study of the
geometry and statistics of the data, its possible decomposition
into clusters, etc \cite{Hastie}. In addition, there are many
problems concerning supervised learning, in which additional
information, such as a discrete class $g(\x)\in\{g_1,...,g_k\}$ or
a continuous function value $f(\x)$ is given to some of the data
points. In this paper we are concerned only with the unsupervised
case, although some of the methods and ideas presented can be
applied to the supervised or semi-supervised case as well
\cite{PNAS1}.

In the later three above-mentioned application fields the data is
typically sampled from a complex biological, chemical or physical
{\em dynamical} system, in which there is an inherent notion of
time. Many of these systems involve multiple time and length
scales, and  in many interesting cases there is a separation of
time scales, that is, there are only a few "slow" time scales at
which the system performs conformational changes from one
meta-stable state to another, with many additional fast time
scales at which the system performs local fluctuations within
these meta-stable states. In the case of macromolecules the slow
time scale is that of a conformational change, while the fast time
scales are governed by the chaotic rotations and vibrations of the
individual chemical bonds between the different atoms of the
molecule, as well as the random fluctuations  due to the frequent
collisions with the surrounding solvent water molecules. In the
more general case of interacting particle systems, the fast time
scales are those of density fluctuations around the mean density
profiles, while the slow time scales correspond to the time
evolution of these mean density profiles.

Although on the fine time and length scales the full description
of such systems requires a high dimensional space, e.g. the
locations (and velocities) of all the different particles, these
systems typically have an intrinsic low dimensionality on coarser
length and time scales. Thus, the coarse time evolution of the
high dimensional system can be described by only a few dynamically
relevant variables, typically called reaction coordinates. Important tasks in such systems are
the reduction of the dimensionality at these coarser scales (known
as homogenization), and the efficient representation of the
complicated linear or non-linear operators that govern their
(coarse grained) time evolution. Additional goals are the
identification of the meta-stable states, the characterization of
the transitions between them and the efficient computation of mean
exit times, potentials of mean force and effective diffusion
coefficients \cite{GKS,Huisinga,Huisinga03,Elber}.

In this paper, following \cite{Lafon}, we consider a family of
diffusion maps for the analysis of these problems. Given a large
dataset, we construct a family of random walk processes based on
isotropic and anisotropic diffusion kernels and study their first
few eigenvalues and eigenvectors (principal components). The key
point in our analysis is that these eigenvectors and eigenvalues
capture important geometrical and statistical information
regarding the structure of the underlying datasets.

It is interesting to note that similar approaches have been
suggested in various different fields. In graph theory, the first
few eigenvectors of the normalized graph Laplacian have been used
for spectral clustering \cite{Weiss99,Weiss}, approximations to
the optimal normalized-cut problem \cite{Malik} and dimensionality
reduction \cite{Belkin,Saerens}, to name just a few. Similar
constructions have also been used for the clustering and
identification of meta-stable states for datasets sampled from
dynamical systems \cite{Huisinga}. However, it seems that the
connection of these computed eigenvectors to the underlying
geometry and probability density of the dataset has not been fully
explored.

In this paper, we consider the connection of these eigenvalues and
eigenvectors to the underlying geometry and probability density
distribution of the dataset. To this end, we assume that the data
is sampled from some (unknown) probability distribution, and view
the eigenvectors computed on the finite dataset as discrete
approximations of corresponding eigenfunctions of suitably defined
continuum operators in an infinite population setting. As the
number of samples goes to infinity, the discrete random walk on
the set converges to a diffusion process defined on the
$n$-dimensional space but with a non-uniform probability density.
By explicitly studying the asymptotic form of the
Chapman-Kolmogorov equations in this setting (e.g., the
infinitesimal generators), we find that for data sampled from a
general probability distribution, written in Boltzmann form as
$p(\x)= e^{-U(\x)}$, the eigenfunctions and eigenvalues of the
standard normalized graph Laplacian construction correspond to a
diffusion process with a potential $2 U(\x)$ (instead of a single
$U(\x)$). Therefore, a subset of the first few eigenfunctions are
indeed well suited for spectral clustering of data that contains
only a few well separated clusters, corresponding to deep wells in
the potential $U(\x)$.

Motivated by the well known connection between diffusion processes
and Schr\"odinger operators \cite{Bernstein}, we propose a
different novel non-isotropic construction of a random walk on the
graph, that in the asymptotic limit of infinite data recovers the
eigenvalues and eigenfunctions of a diffusion process with the
same potential $U(\x)$. This normalization, therefore, is most
suited for the study of the long time behavior of complex
dynamical systems that evolve in time according to a stochastic
differential equation. For example, in the case of a dynamical
system driven by a bistable potential with two wells, (e.g. with
one slow time scale for the transition between the wells and many
fast time scales) the second eigenfunction can serve as a
parametrization of the reaction coordinate between the two states,
much in analogy to its use for the construction of an
approximation to the optimal normalized cut for graph
segmentation. For the analysis of dynamical systems, we also
propose to use a subset of the first few eigenfunctions as
reaction coordinates for the design of fast simulations. The main
idea is that once a parametrization of dynamically meaningful
reaction coordinates is known, and lifting and projection
operators between the original space and the diffusion map are
available, detailed simulations can be initialized at different
locations on the reaction path and run only for short times, to
estimate the transition probabilities to different nearby
locations in the reaction coordinate space, thus efficiently
constructing a potential of mean field and an efficient diffusion
coefficient on the reaction path \cite{Yannis}.

Finally, we describe yet another random walk construction that in
the limit of infinite data recovers the Laplace-Beltrami (heat)
operator on the manifold on which the data resides, regardless of
the possibly non-uniform sampling of points on the manifold. This
normalization is therefore best suited for learning the geometry
of the dataset, as it separates the geometry of the manifold from
the statistics on it.

Our analysis thus reveals the intimate connection between the
eigenvalues and eigenfunctions of different random walks on the
finite graph to the underlying geometry and probability
distribution $p=e^{-U}$ from which the dataset was sampled. These
findings lead to a better understanding of the advantages and
limitations of diffusion maps as a tool to solve different tasks
in the analysis of high dimensional data.

\section{Problem Setup}\label{s:setup}

Consider a finite dataset $\{\x_i\}_{i=1}^N \in \mathbb{R}^n$. We
consider two different possible scenarios for the origin of the
data. In the first scenario, the data is not necessarily derived
from a dynamical system, but rather it is randomly sampled from
some arbitrary probability distribution $p(\x)$. In this case we
define an associated potential
\begin{equation}
  U(\x) = - \log p(\x)
\end{equation}
so that $p= e^{-U}$.

In the second scenario, we assume that the data is sampled from a
dynamical system in equilibrium. We further assume that the
dynamical system, defined by its state $\mb{\x}(t)\in\mathbb{R}^n$
at time $t$, satisfies the following generic stochastic
differential equation (SDE)
\begin{equation}
\dot{\mb{x}} = -\nabla U(\x) + \sqrt{2} \dot{\mb{w}}
    \label{SDE}
\end{equation}
where a dot on a variable means differentiation with respect to
time, $U$ is the free energy at $\x$ (which, with some abuse of
nomenclature, we will also call the potential at $\x$), and
$\mb{w}(t)$ is an $n$-dimensional Brownian motion process. In this
case there is an explicit notion of time, and the transition
probability density $p(\x,t|\y,s)$ of finding the system at
location $\x$ at time $t$, given an initial location $\y$ at time
$s$ ($t > s$), satisfies the forward Fokker-Planck equation (FPE)
\cite{Schuss,Gardiner}
\begin{equation}
  \frac{\p p}{\p t} = \nabla \cdot \left(\nabla p + p \nabla U(\x)\right)
    \label{FPE}
\end{equation}
with initial condition
\begin{equation}
\lim_{t\to s^+}  p(\x,t|\y,s) = \delta(\x - \y)
\end{equation}
Similarly, the backward Fokker-Planck equation for the density
$p(\x,t | \y, s)$, in the backward variables $\y,s$ ($s<t$) is
\begin{equation}
  -\frac{\p p}{\p s } = \Delta p - \nabla p \cdot \nabla U(\y)
    \label{backward_FPE}
\end{equation}
where differentiations in (\ref{backward_FPE}) are with respect to
the variable $\y$, and the Laplacian $\Delta$ is a negative
operator, defined as $\Delta u = \nabla \cdot (\nabla u)$.

As time $t\to\infty$ the steady state solution of (\ref{FPE}) is
given by the equilibrium Boltzmann probability density,
\begin{equation}
\mu(\x)d\x =  \Pr\{\x\}d\x = \frac{\exp(-U(\x))}{Z}d\x
    \label{mu_x}
\end{equation}
where $Z$ is a normalization constant (known as the partition
function in statistical physics), given by
\begin{equation}
  Z = \int_{\mathbb{R}^n} \exp(-U(\x)) d\x 
\end{equation}
In what follows we assume that the potential $U(\x)$ is shifted by
the suitable constant (which does not change the SDE (\ref{SDE})),
so that $Z=1$. Also, we use the notation $\mu(\x) = \Pr\{\x\} =
p(\x)=e^{-U(\x)}$ interchangeably to denote the (invariant)
probability measure on the space.

Note that in both scenarios, the steady state probability density,
given by (\ref{mu_x}) is identical. Therefore, for the purpose of
our initial analysis, which does not directly take into account the
possible time dependence of the data, it is only the features of
the underlying potential $U(\x)$ that come into play.

The Langevin equation (\ref{SDE}) or the corresponding
Fokker-Planck equation (\ref{FPE}) are commonly used to describe
mechanical, physical, chemical, or biological systems driven by
noise. The study of their behavior, and specifically the decay to
equilibrium has been the subject of much theoretical research
\cite{Risken}. In general, the solution of the Fokker-Planck
equation (\ref{FPE}) can be written in terms of an eigenfunction
expansion
\begin{equation}
p(\x,t) = \sum_{j=0}^\infty a_j e^{-\lambda_j t} \varphi_j(\x)
\end{equation}
where $-\lambda_j$ are the eigenvalues of the FP operator, with
$\lambda_0 = 0 < \lambda_1\leq \lambda_2\leq\ldots$,
$\varphi_j(\x)$ are their corresponding eigenfunctions, and the
coefficients $a_j$ depend on the initial conditions. Obviously,
the long term behavior of the system is governed only by the first
few eigenfunctions $\varphi_0,\varphi_1,\ldots,\varphi_k$, where
$k$ is typically small and depends on the time scale of interest.
In low dimensions, e.g. $n\leq 3$ for example, it is possible to
calculate approximations to these eigenfunctions via numerical
solutions of the relevant partial differential equations. In high
dimensions, however, this approach is in general infeasible and
one typically resorts to simulations of trajectories of the
corresponding SDE (\ref{SDE}). In this case, there is a need to
employ statistical methods to analyze the simulated trajectories,
identify the slow variables, the meta-stable states, the reaction
pathways connecting them and the mean transition times between
them.

\section{Diffusion Maps}

\subsection{Finite Data}\label{s:discrete}

Let $\{\x_i\}_{i=1}^N$, denote the $N$ samples, either merged from
many different simulations of the stochastic equation (\ref{SDE}),
or simply given without an underlying dynamical system. In
\cite{Lafon}, Coifman and Lafon suggested the following method,
based on the definition of a Markov chain on the data, for the
analysis of the geometry of general datasets:

For a fixed value of $\ve$ (a metaparameter of the algorithm),
define an isotropic diffusion kernel,
\begin{equation}
  k_\ve(\x,\y) =
  \exp\left(-\frac{\norm{\x-\y}^2}{2\ve}\right) \label{k_epsilon}
\end{equation}
Assume that the transition probability between points $\x_i$ and
$\x_j$ is proportional to $k_\ve(\x_i,\x_j)$, and construct an
$N\times N$ Markov matrix, as follows
\begin{equation}
  M(i,j) = \frac{k_\ve(\x_i,\x_j)}{p_\ve(\x_j)}
    \label{M_discrete}
\end{equation}
where $p_\ve$ is the required normalization constant, given by
\begin{equation}
  p_\ve(\x_j) = \sum_i k_\ve(\x_i,\x_j)
    \label{p_ve_discrete}
\end{equation}
For large enough values of $\ve$ the Markov matrix $M$ is fully
connected (in the numerical sense) and therefore has an eigenvalue
$\lambda_0=1$ with multiplicity one and a sequence of additional
$n-1$ non-increasing eigenvalues $\lambda_j < 1$, with
corresponding eigenvectors $\varphi_j$.

The diffusion map at time $m$ is defined as the mapping from $\x$
to the vector
\[
\Phi_m(\x) = \left(\lambda_0^m
\varphi_0(\x),\lambda_1^m\varphi_1(\x),\ldots,\lambda_k^m
\varphi_k(\x)\right) \] for some small value of $k$. In
\cite{Lafon}, it was demonstrated that this mapping gives a low
dimensional parametrization of the geometry and density of the
data. In the field of data analysis, this construction is known as
the {\em normalized graph Laplacian}. In \cite{Malik}, Shi and
Malik suggested using the first non-trivial eigenvector to compute
an approximation to the optimal normalized cut of a graph, while
the first few eigenvectors were suggested by Weiss et al.
\cite{Weiss99,Weiss} for clustering. Similar constructions,
falling under the general term of kernel methods have been used in
the machine learning community for classification and regression
\cite{Kernel}. In this paper we elucidate the connection between
this construction and the underlying potential $U(\x)$.

\subsection{The Continuum Diffusion Process}\label{sec: continuum diffusion process}

To analyze the eigenvalues and eigenvectors of the normalized
graph Laplacian, we consider them as a finite approximation of a
suitably defined diffusion operator acting on the continuum
probability space from which the data was sampled. We thus
consider the limit of the above Markov chain process as the number
of samples approaches infinity. For a finite value of $\ve$, the
Markov chain in discrete time and space converges to a Markov
process in discrete time but continuous space. Then, in the limit
$\ve\to0$, this jump process converges to a diffusion process on
$\mathbb{R}^n$, whose local transition probability depends on the
non-uniform probability measure $\mu(\x) = e^{-U(\x)}$.

We first consider the case of a fixed $\ve > 0$, and take $N\to\infty$. Using the
similarity of (\ref{k_epsilon}) to the diffusion kernel, we view
$\ve$ as a measure of time and consider a discrete jump process at
time intervals $\Delta t= \ve$, with a transition probability
between points $\y$ and $\x$ proportional to $k_\ve(\x,\y)$.
However, since the density of points is not uniform but rather
given by the measure $\mu(\x)$, we define an associated
normalization factor $p_\ve(\y)$ as follows,
\begin{equation}
p_\ve(\y) = \int k_\ve(\x,\y) \mu(\x) d\x
    \label{p_ve}
\end{equation}
and a forward transition probability
\begin{equation}
M_f(\x|\y) = \Pr(\x(t+\ve) = \x\,|\x(t)=\y) =
\frac{k_\ve(\x,\y)}{p_\ve(\y)}
    \label{M_f}
\end{equation}
Equations (\ref{p_ve}) and (\ref{M_f}) are the continuous
analogues of the discrete equations (\ref{p_ve_discrete}) and
(\ref{M_discrete}). For future use, we also define a symmetric
kernel $M_s(\x,\y)$ as follows,
\begin{equation}
  M_s(\x,\y) = \frac{k_\ve(\x,\y)}{\sqrt{p_\ve(\x)p_\ve(\y)}}
    \label{M_s}
\end{equation}
Note that $p_\ve(\x)$ is an estimate of the local probability
density at $\x$, computed by averaging the density in a
neighborhood of radius $O(\ve^{1/2})$ around $\x$. Indeed, as
$\ve\to 0$, we have that
\begin{equation}
  p_\ve(\x) = p(\x) + \frac\ve{2} \Delta p(\x) + O(\ve^{3/2})
\end{equation}

 We now define forward, backward and symmetric Chapman-Kolmogorov operators
on functions defined on this probability space, as follows,
\begin{equation}
  T_f[\varphi](\x) = \int M_f(\x|\y) \varphi(\y) d\mu(\y)
\end{equation}
\begin{equation}
  T_b[\varphi](\x) = \int M_f(\y|\x) \varphi(\y) d\mu(\y)
\end{equation}
and
\begin{equation}
  T_s[\varphi](\x) = \int M_s(\x,\y) \varphi(\y) d\mu(\y)
\end{equation}
If $\varphi(\x)$ is the probability of finding the system at
location $\x$ at time $t=0$, then $T_f[\varphi]$ is the evolution
of this probability to time $t=\ve$. Similarly, if $\psi(\mb{z})$
is some function on the space, then $T_b[\psi](\x)$ is the mean
(average) value of that function at time $\ve$ for a random walk
that started at $\x$, and so $T_b^m[\psi](\x)$ is the average
value of the function at time $t=m\ve$.

By definition, the operators $T_f$ and $T_b$ are adjoint under the
inner product with weight $\mu$, while the operator $T_s$ is self
adjoint under this inner product,
\begin{equation}
  \langle T_f \varphi , \psi \rangle_{\mu} = \langle \varphi, T_b \psi
  \rangle_{\mu},
        \quad \quad
  \langle T_s \varphi , \psi \rangle_{\mu} = \langle \varphi, T_s \psi
  \rangle_{\mu}
\end{equation}
Moreover, since $T_s$ is obtained via conjugation of the kernel
$M_f$ with $\sqrt{p_\ve(\x)}$ all three operators share the same
eigenvalues. The corresponding eigenfunctions can be found via
conjugation by $\sqrt{p_\ve}$. For example, if $T_s\varphi_s =
\lambda \varphi_s$, then the corresponding eigenfunctions for
$T_f$ and $T_b$ are $\varphi_f = \sqrt{p_\ve} \varphi_s$ and
$\varphi_b = \varphi_s/\sqrt{p_\ve}$, respectively. Since
$\sqrt{p_\ve}$ is the first eigenfunction with $\lambda_0 = 1$ of
$T_s$, the steady state of the forward operator is simply
$p_\ve(\x)$, while for the backward operator it is the uniform
density, $\varphi_b=1$.

Obviously, the eigenvalues and eigenvectors of the discrete Markov
chain described in the previous section are discrete
approximations to the eigenvalues and eigenfunctions of these
continuum operators. Rigorous mathematical proofs of this
convergence as $N\to\infty$ under various assumptions have been
recently obtained by several authors \cite{BelkinC,Hein}.
Therefore, for a better understanding of the finite sample case,
we are interested in the properties of the eigenvalues and
eigenfunctions of either one of the operators $T_f,T_b$ or $T_s$,
and how these relate to the measure $\mu(\x)$ and to the
corresponding potential $U(\x)$. To this end, we look for
functions $\varphi(\x)$ such that
\begin{equation}
T_j\varphi =  \int M_j(\x,\y) \varphi(\y) \Pr\{\y\} d\y = \lambda
\varphi(\x)
    \label{Tlambda}
\end{equation}
where $j \in\{f,b,s\}$.

While in the case of a finite amount of data, $\ve$ must remain
finite so as to have enough neighbors in a ball of radius
$O(\ve^{1/2})$ near each point $\x$, as the number of samples goes
to infinity, we can take smaller and smaller values of $\ve$.
Therefore, it is instructive to look at the limit $\ve \to 0$. In
this case, the transition probability densities of the continuous
in space but discrete in time Markov chain converge to those of a
diffusion process, whose time evolution is described by a
differential equation
\[
  \frac{\p \varphi}{\p t} = {\cal H}_f \varphi
\]
where ${\cal H}_f$ is the infinitesimal generator or propagator of
the forward operator, defined as
\[
{\cal  H}_f = \lim_{\ve \to 0}\frac{I - T_f}{\varepsilon}
\]
As shown in the Appendix, by computing the asymptotic expansion of
the corresponding integrals in the limit $\ve\to0$, we obtain that
\begin{equation}
{\cal H}_f \varphi = \Delta \varphi
   - \varphi \left(e^U\Delta e^{-U}\right)
    \label{H_f}
\end{equation}
Similarly, the inifinitesimal operator of the backward operator is
given by
\begin{equation}
{\cal  H}_b \psi = \lim_{\ve \to 0} \frac{T_b-I}{\ve}\psi = \Delta
\psi - 2
  \nabla \psi \cdot \nabla U
    \label{H_b}
\end{equation}
As expected, $\psi_0=1$ is the eigenfunction with $\lambda_0=0$ of
the backward infinitesimal operator, while $\varphi_0=e^{-U}$ is
the eigenfunction of the forward one.

A few important remarks are due at this point. First, note that
the backward operator (\ref{H_b}) has the same functional form as
the backward FPE (\ref{backward_FPE}), but with a potential $2
U(\x)$ instead of $U(\x)$.  The forward operator (\ref{H_f}) has a
different functional form from the forward FPE (\ref{FPE})
corresponding to the stochastic differential equation (\ref{SDE}).
This should come as no surprise, since (\ref{H_f}) is the
differential operator of an isotropic diffusion process on a space
with non-uniform probability measure $\mu(\x)$, which is obviously
different from the standard anisotropic diffusion in a space with
a uniform measure, as described by the SDE (\ref{SDE})
\cite{Gardiner}.

Interestingly, however, the form of the forward operator is
the same as the Schr\"{o}dinger operator of quantum physics
\cite{Singh}, e.g.
\begin{equation}
  {\cal H}\varphi = \Delta \varphi - \varphi V(\x) \label{schrodinger}
    \label{QM}
\end{equation}
where in our case the potential $V(\x)$ has the following specific
form
\begin{equation}
  V(\x) =\left(\nabla U(\x)\right)^2 - \Delta U(\x).
    \label{Vx}
\end{equation}
Therefore, in the limit $N  \to \infty, \ve\to 0$, the
eigenfunctions of the diffusion map are the same as those of the
Schr\"odinger operator (\ref{schrodinger}) with a potential
(\ref{Vx}). The properties of the first few of these
eigenfunctions have been extensively studied theoretically for
simple potentials $V(\x)$ \cite{Singh}.

In order to see why the forward operator ${\cal H}_f$ also
corresponds to a potential $2U(\x)$ instead of $U(\x)$, we recall
that there is a well known correspondence \cite{Bernstein},
between the Schr\"{o}dinger equation with a sypersymmetric
potential of the specific form (\ref{Vx}) and a diffusion process
described by a Fokker-Planck equation of the standard form
(\ref{FPE}). The transformation
\begin{equation}
  p(\x,t) = \psi(\x,t) e^{-U(\x)/2}
    \label{transformation}
\end{equation}
applied to the original FPE (\ref{FPE}) yields the Schr\"odinger
equation with imaginary time
\begin{equation}
-  \frac{\p \psi}{\p t}  = \Delta \psi - \psi\left(\frac{(\nabla
U)^2}4 - \frac{\Delta U}2\right)\label{eq:Schrodinger imaginary time}
\end{equation}
    Comparing (\ref{eq:Schrodinger imaginary time}) with (\ref{Vx}),
we conclude that the eigenvalues of the operator (\ref{H_f}) are
the same as those of a Fokker-Planck equation with a potential $2
U(\x)$. Therefore, in general, for data sampled from the SDE
(\ref{SDE}), there is no direct correspondence between the
eigenvalues and eigenfunctions of the normalized graph Laplacian
and those of the corresponding Fokker-Planck equation (\ref{FPE}).
However, when the original potential $U(\x)$ has two metastable
states separated by a large barrier, corresponding to two well
separated clusters, so does $2U(\x)$. Therefore, the first
non-trivial eigenvalue is governed by the mean passage time
between the two barriers, and the first non-trivial eigenfunction
gives a parametrization of the path between them (see also the
analysis of the next section).

We note that in \cite{Horn}, Horn and Gottlieb suggested a
clustering algorithm based on the Schr\"{o}dinger operator
(\ref{QM}), where they constructed an approximate eigenfunction
$\psi(\x) = p_\ve(\x)$ as in our eq. \ref{p_ve_discrete}), and
computed its corresponding potential $V(\x)$ from eq. (\ref{QM}).
The clusters were then defined by the minima of the potential $V$.
Employing the same asymptotic analysis of this paper, one can show
that in the appropriate limit, the computed potential $V$ is given
by (\ref{Vx}). This asymptotic analysis and the connection between
the quantum operator and a diffusion process thus provides a
mathematical explanation for the success of their method. Indeed,
when $U$ has a deep parabolic minima at a point $\x$,
corresponding to a well defined cluster, so does $V$.

\section{Anisotropic Diffusion Maps}\label{ref: anisotropic diffusion maps}

In the previous section we showed that asymptotically, the
eigenvalues and eigenfunctions of the normalized graph Laplacian
operator correspond to the Fokker-Planck equation with a potential
$2U(\x)$ instead of the single $U(\x)$. In this section we present
a different normalization that yields infinitesimal generators
corresponding to the potential $U(\x)$ without the additional
factor of two.

In fact, following \cite{Lafon}
we consider in more generality a whole family of
different normalizations and their corresponding diffusions, and
we show that, in addition to containing the graph Laplacian
normalization used in the previous section, this collection of
diffusions includes at least two other Laplacians of interest: the
Laplace-Beltrami operator, which captures the Riemannian geometry
of the data set, and the backward Fokker-Planck operator of
equation (\ref{backward_FPE}).

Instead of applying the graph Laplacian normalization to the
isotropic kernel $k_\varepsilon(\x,\y)$, we first appropriately
renormalize the kernel into an anisotropic one to obtain a new
weighted graph, and then apply the graph Laplacian normalization
to this graph. More precisely, we proceed as follows: start with a
Gaussian kernel $k_\varepsilon(\x,\y)$ and let $\alpha>0$ be a
parameter indexing our family of diffusions. Define an estimate
for the local density as
\[
p_\ve(\x)=\int k_\varepsilon(\x,\y) \Pr\{\y\}d\y
\]
and consider the family of kernels
\[
k^{(\alpha)}_\varepsilon(\x,\y)=\frac{k_\varepsilon(\x,\y)}{p_\varepsilon^\alpha(\x)p_\varepsilon^\alpha(\y)}
\]
We now apply the graph Laplacian normalization by computing the
normalization factor
\[
d_\ve^{(\alpha)}(\y)=\int
k^{(\alpha)}_\varepsilon(\x,\y)\Pr\{\x\}d\x
\]
and forming a forward transition probability kernel
\[
M_f^{(\alpha)}(\x|\y)=\Pr\{\x(t+\varepsilon)=\x|\x(t)=\y\}=\frac{k_\ve^{(\alpha)}(\x,\y)}{d_\ve^{(\alpha)}(\y)}
\]
Similar to the analysis of section \ref{sec: continuum diffusion
process}, we can construct the corresponding forward, symmetric
and backward diffusion kernels. It can be shown (see appendix
\ref{infinitesimal computations}) that the forward and backward
infinitesimal generators of this diffusion are
\begin{eqnarray}
\mathcal H_b^{(\alpha)}\psi &=& \Delta \psi -
2(1-\alpha)\nabla\phi\cdot \nabla U \\
 \mathcal H_f^{(\alpha)} \varphi&=&\Delta
\varphi-2\alpha \nabla \varphi \cdot \nabla U + (2\alpha-1)
\varphi \left((\nabla U)^2 - \Delta U\right)
\end{eqnarray}
We mention three interesting cases:
\begin{itemize}
\item For $\alpha=0$, this construction yields the classical
normalized graph Laplacian with the infinitesimal operator given
by equation (\ref{H_f})
\[
\mathcal H_f \varphi=\Delta \varphi-\left(e^{U}\Delta
e^{-U}\right)\varphi
\]
\item For $\alpha=1$, the backward generator gives the
Laplace-Beltrami operator:
\begin{equation}
\mathcal H_b\psi=\Delta \psi
\end{equation}
In other words, this diffusion captures only the geometry of the
data, in which the density $e^{-U}$ plays absolutely no role.
Therefore, this normalization separates the geometry of the
underlying manifold from the statistics on it.

\item For $\alpha=\frac 1 2$, the infinitesimal operator of the
forward and backward operators coincide and are given by
\begin{equation}
\mathcal H_f \varphi= \mathcal H_b \varphi = \Delta \varphi-
\nabla \varphi \cdot \nabla U
\end{equation}
which is exactly the backward FPE (\ref{backward_FPE}), with a
potential $U(\x)$.
\end{itemize}

Therefore, the last case with $\alpha=1/2$ provides a consistent
method to approximate the eigenvalues and eigenfunctions
corresponding to the stochastic differential equation (\ref{SDE}).
This is done by constructing a graph Laplacian with an
appropriately anisotropic weighted graph.

As explained in \cite{Lafon,Saerens,new}, the Euclidian distance
between any two points after the diffusion map embedding into
$\mathbb{R}^k$ is almost equal to their diffusion distance on the
original dataset. Thus, for dynamical systems with only one or two
slow time scales, and many fast time scales, only a small number
of diffusion map coordinates need be retained for the coarse
grained representation of the data at medium to long times, at
which the fast coordinates have equilibrated. Therefore, the
diffusion map can be considered as an empirical method to perform
data-driven or equation-free homogenization. In particular, since
this observation yields a computational method for the
approximation of the top eigenfunctions and eigenvalues, this
method can be applied towards the design of fast and efficient
simulations that can be initialized on arbitrary points on the
diffusion map. This application will be described in a separate
publication \cite{new}.

\section{Examples}

In this section we present the potential strength of the diffusion
map method by analyzing, both analytically and numerically a few
toy examples, with simple potentials $U(\x)$. More complicated
high dimensional examples of stochastic dynamical systems are
analyzed in \cite{new}, while other applications such as the
analysis of images for which we typically have no underlying
probability model appear in \cite{Lafon}.

\subsection{Parabolic potential in 1-D}

We start with the simplest case of a parabolic potential in one
dimension, which in the context of the SDE (\ref{SDE}),
corresponds to the well known Ornstein-Uhlenbeck process. We thus
consider a potential $U(x) = x^2 /2 \tau$, with a corresponding
normalized density $p = e^{-U}/\sqrt{2\pi\tau}$.

The normalization factor $p_\ve$ can be computed explicitly
\[
  p_\ve(y) = \int \frac{e^{-(x-y)^2/2\ve}}{\sqrt{2\pi\ve}}
  \frac{e^{-x^2/2\tau}}{\sqrt{2\pi\tau}}dx = \frac{1}{\sqrt{2\pi(\tau + \ve)}}
  e^{-y^2/2(\tau+\ve)}
\]
where, for convenience, we multiplied the kernel $k_\ve(x,y)$ by a
normalization factor $1/\sqrt{2\pi\ve}$. Therefore, the
eigenvalue/eigenfunction problem for the symmetric operator $T_s$
with a finite $\ve$ reads
\[
T_s\varphi =   \int
\frac{\exp\left(-\frac{(x-y)^2}{2\varepsilon}\right)}{\sqrt{2\pi\ve}}
  \exp\left(\frac{x^2+y^2}{4(\varepsilon+\tau)}\right)
\exp\left(-\frac{y^2}{2\tau}\right) \sqrt{\frac{\tau+\ve}{\tau}}
  \varphi(y)dy = \lambda \varphi(x)
\]
The first eigenfunction, with eigenvalue $\lambda_0=1$ is given by
\[
  \varphi_0(x) = C \sqrt{p_\ve(x)} = C \exp\left(-\frac{x^2}{4(\varepsilon + \tau)}\right)
\]
The second eigenfunction, with eigenvalue $\lambda_1 = \tau/(\tau
+ \varepsilon) < 1$ is, up to normalization
\[
  \varphi_1(x) = x \exp\left(-\frac{x^2}{4(\varepsilon + \tau)}\right)
\]
In general, the sequence of eigenfunctions and eigenvalues is
characterized by the following lemma:

\noindent {\bf \em Lemma:} The eigenvalues of the operator $T_s$
are $\lambda_k = \left(\tau/(\tau+\varepsilon)\right)^k$, with the
corresponding eigenvectors given by
\begin{equation}
\varphi_k(x) = p_k(x)
\exp\left(-\frac{x^2}{4(\tau+\varepsilon)}\right)
\end{equation}
where $p_k$ is a polynomial of degree $k$ (even or odd depending
on $k$).

In the limit $\ve\to 0$, we obtain the eigenfunctions of the
corresponding infinitesimal generator. For the specific potential
$U(x)=x^2/2\tau$, the eigenfunction problem for the backward
generator reads
\begin{equation}
  \psi_{xx} - 2 \frac{x}{\tau} \psi_x = - \lambda \psi
\end{equation}
The solutions of this eigenfunction problem are, up to scaling of
$x$, the well known Hermite polynomials, which by the
correspondence of this operator to the Schr\"{o}dinger
eigenvector/eigenvalue problem, are also the eigenfunctions of the
quantum harmonic oscillator (after multiplication by the
appropriate Gaussian) \cite{Singh}.

Note that plotting the second vs. the first eigenfunctions (with
the convention that the zeroth eigenfunction is the constant one,
which we typically ignore), is the same as plotting $x^2+1$ vs
$x$, e.g. a parabola. Therefore, we expect that for a large enough
and yet finite data-set sampled from this potential, the plot of
the corresponding discrete eigenfunctions should lay on a
parabolic curve (see next section for a numerical example).

\subsection{Multi-Dimensional Parabolic Potential}

We now consider a harmonic potential in $n$-dimensions, of the
form
\begin{equation}
  U(\x) = \sum_j \frac{x_j^2}{2\tau_j}
\end{equation}
where, in addition, we assume $\tau_1\gg \tau_2 > \tau_3 >\ldots >
\tau_n$, so that $x_1$ is a slow variable in the context of the
SDE (\ref{SDE}).

We note that for this specific potential, the probability density
has a separable  structure, $p(\x) = p_1(x_1)\ldots p_n(x_n)$, and
so does the kernel $k_\ve(\x,\y)$, and consequently, also the
symmetric kernel $M_s(\x,\y)$. Therefore, there is an
outer-product structure to the eigenvalues and eigenfunctions. For
example, in two dimensions the eigenfunctions and eigenvalues are
given by
\begin{equation}
\varphi_{i,j}(x_1,x_2) = \varphi_{1,i}(x_1)\varphi_{2,j}(x_2)\quad
\mbox{and}\quad\lambda_{i,j} = \mu_1^i \mu_2^j
\end{equation}
where $\mu_1 = \tau_1/(\tau_1+\ve)$ and $\mu_2 = \tau_2/(\tau_2
+\ve)$. Since by assumption $\tau_1 \gg \tau_2$, then upon
ordering of the eigenfunctions by decreasing eigenvalue, the first
non-trivial eigenfunctions are
$\varphi_{1,0},\varphi_{2,0},\ldots$, which depend only on the
slow variable $x_1$. Note that indeed, regardless of the value of
$\ve$, as long as $\tau_2 > 2 \tau_1$, we have that $\lambda_1^2 >
\lambda_2$. Therefore, in this example the first few coordinates
of the diffusion map give a (redundant) parametrization of the
slow variable $x_1$ in the system.

\begin{figure}[t]
\mbox{
\begin{minipage}[t] {\textwidth}
\begin{center}
\begin{tabular}{c}
\psfig{figure=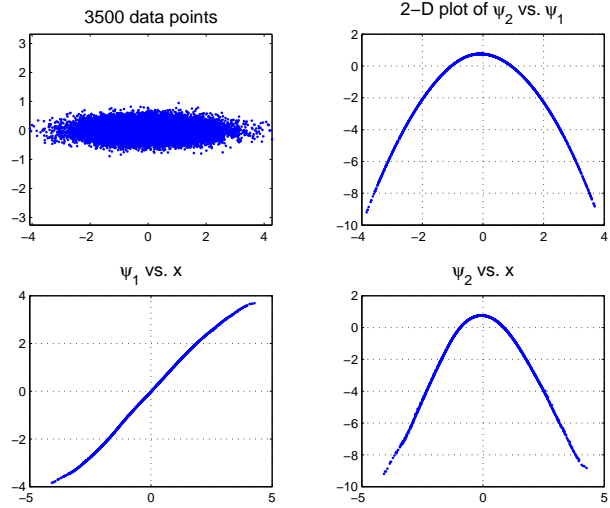,width=8.0cm}\\
\end{tabular}
\end{center}
\end{minipage}
}
\\
\caption{The anisotropic diffusion map on a harmonic potential in
2-D. } \label{f:u1}
\end{figure}

In figure \ref{f:u1} we present numerical results corresponding to
a 2-dimensional potential with $\tau_1=1,\tau_2=1/25$. In the
upper left some 3500 points sampled from the distribution
$p=e^{-U}$ are shown. In the lower right and left panels, the
first two non-trivial backward eigenfunctions $\psi_1$ and
$\psi_2$ are plotted vs. the slow variable $x_1$. Note that except
at the edges, where the statistical sampling is poor, the first
eigenfunction is linear in $x_1$, while the second one is
quadratic in $x_1$. In the upper right panel we plot $\psi_2$ vs.
$\psi_1$ and note that they indeed lie on a parabolic curve, as
predicted by the analysis of the previous section.

\subsection{A potential with two minima}

\begin{figure}[t]
\mbox{
\begin{minipage}[t] {\textwidth}
\begin{center}
\begin{tabular}{c}
\psfig{figure=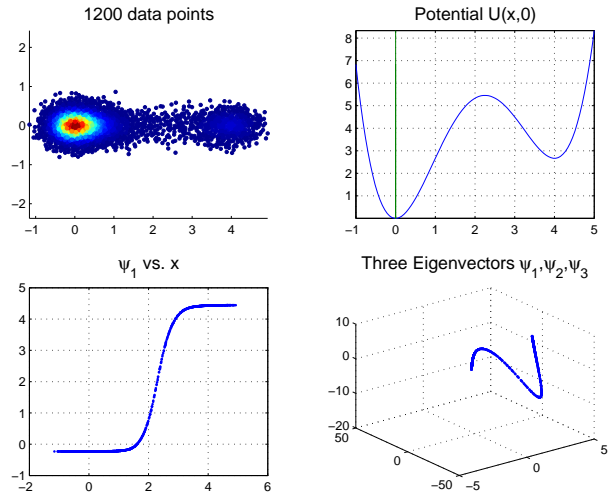,width=8.0cm}\\
\end{tabular}
\end{center}
\end{minipage}
}
\\
\caption{Numerical results for a double well potential in 2-D. }
\label{f:u2}
\end{figure}

We now consider a double well potential $U(x)$ with two minima,
one at $x_L$ and one at $x_R$. For simplicity of the analysis, we
assume a symmetric potential around $(x_L+x_R)/2$, with
$U(x_L)=U(x_R) = 0$ (see figure \ref{f:u2}). In the context of
data clustering, this can be viewed as approximately a mixture of
two Gaussian clouds, while in the context of stochastic dynamical
systems, this potential defines two meta-stable states.

We first consider an approximation to the quantity $p_\ve(x)$,
given by eq. (\ref{p_ve}). For $x$ near $x_L$, $U(x) \approx
(x-x_L)^2/\tau_L$, while for $x$ near $x_R$, $U(x) \approx
(x-x_R)^2/\tau_R$. Therefore,
\begin{equation}
  e^{-U(y)} \approx e^{-(y-x_L)^2/2\tau_L} + e^{-(y-x_R)^2/2\tau_R}
\end{equation}
and
\begin{eqnarray}
  p_\ve(x) &\approx& \frac{1}{\sqrt{2}} \left(\frac{\sqrt{\tau_L}}{\sqrt{\tau_L+\ve}}e^{-(x-x_L)^2/2(\tau_L + \ve)}
+\frac{\sqrt{\tau_R}}{\sqrt{\tau_R+\ve}}e^{-(x-x_R)^2/2(\tau_R+
\ve)}
  \right) \nonumber \\
    &=& \frac1{\sqrt{2}} \left[\varphi_L(x) + \varphi_R(x)\right]
\end{eqnarray}
where $\varphi_L$ and $\varphi_R$ are the first forward
eigenfunctions corresponding to a single well potential centered
at $x_L$ or at $x_R$, respectively. As is well known both in the
theory of quantum physics and in the theory of the Fokker-Planck
equation, an approximate expression for the next eigenfunction is
\[
\varphi_1(x) = \frac1{\sqrt{2}} \left[\varphi_L(x) -
\varphi_R(x)\right]
\]
Therefore, the first non-trivial eigenfunction of the backward
operator is given by
\[
\psi_1(x) = \frac{\varphi_L(x) - \varphi_R(x)}{\varphi_L(x) +
\varphi_R(x)}
\]
This eigenfunction is roughly $+1$ in one well and $-1$ in the
other well, with a sharp transition between the two values near
the barrier between the two wells. Therefore, this eigenfunction
is indeed suited for clustering. Moreover, in the context of a
mixture of two Gaussian clouds, clustering according to the sign
of $\psi_1(x)$ is asymptotically equivalent to the optimal Bayes
classifier.

\noindent {\bf Example:} Consider the following potential in two
dimensions,
\begin{equation}
  U(x,y) = \frac1{4}\,x^4-\frac{25}{12}x^3+\frac9{2}x^2 + 25
  \frac{y^2}2
\end{equation}
In the $x$ direction, this potential has a double well shape with
two minima, one at $x=0$ and one at $x=4$, separated by a
potential barrier with a maximum at $x=2.25$.

In figure \ref{f:u2} we show some numerical results of the
diffusion map on some 1200 points sub-sampled from a stochastic
simulation with this potential which generated about 40,000
points. On the upper right panel we see the potential $U(x,0)$,
showing the two wells. In the upper left, a scatter plot of all
the points, color coded according to the value of the local
estimated density $p_\ve$, (with $\ve=0.25$) is shown, where the
two clusters are easily observed. In the lower left panel, the
first non-trivial eigenfunction is plotted vs. the first
coordinate $x$. Note that even though there is quite a bit of
variation in the $y$-variable inside each of the wells, the first
eigenfunction $\psi_1$ is essentially a function of only $x$,
regardless of the value of $y$. In the lower right we plot the
first three backward eigenfunctions. Note that they all lie on a
curve, indicating that the long time asymptotics are governed by
the passage time between the two wells and not by the local
fluctuations inside them.

\subsection{Potential with three wells}

\begin{figure}[t]
\mbox{
\begin{minipage}[t] {\textwidth}
\begin{center}
\begin{tabular}{c}
\psfig{figure=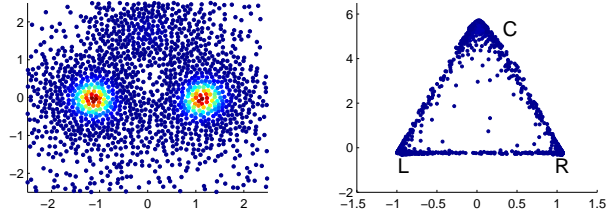,width=8.0cm}\\
\end{tabular}
\end{center}
\end{minipage}
}
\\
\caption{Numerical results for a triple well potential in 2-D. }
\label{f:three_wells}
\end{figure}

We now consider the following two dimensional potential energy
with three wells,
\begin{equation}
  U(x,y) = 3\beta e^{-x^2}\left[e^{-(y-1/3)^2} - e^{-(y-5/3)^2}\right]
  -5\beta
  e^{-y^2}
  \left[e^{-(x-1)^2} + e^{-(x+1)^2}\right]
\end{equation}
where $\beta=1/kT$ is a thermal factor. This potential has two
deep wells at $(-1,0)$ and at $(1,0)$ and a shallower well at
$(0,5/3)$, which we denote as the points $L,R,C$, respectively,
The transitions between the wells of this potential have been
analyzed in many works \cite{Schulten}. In figure
\ref{f:three_wells} we plotted on the left the results of 1400
points sub-sampled from a total of 80000 points randomly generated
from this potential confined to the region $[-2.5,2.5]^2\subset
\mathbb{R}^2$ at temperature $\beta=2$, color-coded by their local
density. On the right we plotted the first two diffusion map
coordinates $\psi_1(\x),\psi_2(\x)$. Notice how in the diffusion
map space one can clearly see a triangle where each vertex
corresponds to one of the points $L,R,C$. This figure shows very
clearly that there are two possible pathways to go from $L$ to
$R$. A direct (short) way and an indirect longer way, that passes
through the shallow well centered at $C$.

\subsection{Iris data set}

\begin{figure}[t]
\mbox{
\begin{minipage}[t] {\textwidth}
\begin{center}
\begin{tabular}{c}
\psfig{figure=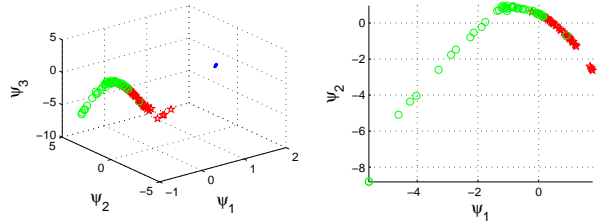,width=8.0cm}\\
\end{tabular}
\end{center}
\end{minipage}
}
\\
\caption{Diffusion map for the iris data set. } \label{f:iris}
\end{figure}

We conclude this section with a diffusion map analysis of one of
the most popular multivariate datasets in pattern recognition, the
iris data set. This set contains 3 distinct classes of samples in
four dimensions, with 50 samples in each class. In figure
\ref{f:iris} we see on the left the result of the three
dimensional diffusion map on this dataset. This picture clearly
shows that all 50 points of class 1 (blue) are shrunk into a
single point in the diffusion map space and can thus be easily
distinguished from classes two and three (red and green). In the
right plot we see the results of re-running the diffusion map on
the 100 remaining red and green samples. The 2-D plot of the first
two diffusion maps coordinates shows that there is no perfect
separation between these two classes. However, clustering
according to the sign of $\psi_1(\x)$ gives misclassifications
rates similar to those of other methods, of the order of 6-8
samples depending on the value chosen for the kernel width $\ve$.

\section{Summary and Discussion}

In this paper, we introduced a mathematical framework for the
analysis of diffusion maps, via their corresponding infinitesimal
generators. Our results show that diffusion maps are a natural
method for the analysis of the geometry and probability
distribution of empirical data sets. The identification of the
eigenvectors of the Markov chain as discrete approximations to the
corresponding differential operators provides a mathematical
justification for their use as a dimensional reduction tool and
gives an alternative explanation for their empirical success in
various data analysis applications, such as spectral clustering
and approximations of optimal normalized cuts on discrete graphs.

We generalized the standard construction of the normalized graph
Laplacian to a one-parameter family of graph Laplacians that
provides a low-dimensional description of the data combining the
geometry of the set with the probability distribution of the data
points. The choice of the diffusion map depends on the task at
hand. If, for example, data points are known to approximately lie
on a manifold, and one is solely interested in recovering the
geometry of this set, then an appropriate normalization of a
Gaussian kernel allows to approximate the Laplace-Beltrami
operator, regardless of the density of the data points. This
construction achieves a complete separation of the underlying
geometry, represented by the knowledge of the Laplace operator,
from the statistics of the points. This is important in situations
where the density is meaningless, and yet points on the manifold
are not sampled uniformly on it. In a different scenario, if the
data points are known to be sampled from the equilibrium
distribution of a Fokker-Planck equation, the long-time dynamics
of the density of points can be recovered from an appropriately
normalized random walk process. In this case, there is a subtle
interaction between the distribution of the points and the
geometry of the data set, and one must correctly account for the
density of the points.

While in this paper we analyzed only Gaussian kernels, our
asymptotic results are valid for general kernels, with the
appropriate modification that take into account the mean and
covariance matrix of the kernel. Note, however, that although
asymptotically in the limit $N\to\infty$ and $\ve\to 0$, the
choice of the isotropic kernel is unimportant, for a finite data set the
choice of both $\ve$ and the kernel can be crucial for the success
of the method.

Finally, in the context of dynamical systems, we showed that
diffusion maps with the appropriate normalization constitute a
powerful tool for the analysis of systems exhibiting different
time scales. In particular, as shown in the different examples,
these time scales can be separated and the long time dynamics can
be characterized by the top eigenfunctions of the diffusion
operator. Last, our analysis paves the way for fast simulations of
physical systems by allowing larger integration steps along slow
variable directions. The exact details required for the design of
fast and efficient simulations based on diffusion maps will be
described in a separate publication \cite{new}.

\noindent{\bf Acknowledgments:} The authors would like to thank
the referee for helpful suggestions and for pointing out ref.
\cite{Horn}.

\appendix

\section{Infinitesimal operators for a family of graph Laplacians}\label{infinitesimal computations}

In this appendix, we present the calculation of the infinitesimal
generators for the different diffusion maps characterized by a
parameter $\alpha$.

Suppose that the data set $X$ consists of a Riemannian manifold
with a density $p(\x)=e^{-U(\x)}$ and let $k_\ve(\x,\y)$ be a
Gaussian kernel. It was shown in \cite{Lafon} that if
$k_\varepsilon$ is scaled appropriately, then for any function
$\phi$ on $X$,
\[ \int_X k_\varepsilon
(\x,\y)\phi(\y)dy=\phi(\x)+\varepsilon (\Delta
\phi(x)+q(\x)\phi(\x)) +O(\varepsilon^{\frac 3 2})
\]
where $q$ is a function that depends on the Riemannian geometry of
the manifold and its embedding in $\mathbb{R}^n$. Using the notations introduced in section \ref{ref:
anisotropic diffusion maps}, it is easy to verify that
\[
p_\ve(\x)=p(\x)+\varepsilon(\Delta p(\x)+q(\x)p(\x)) +
O(\ve^{3/2})
\]
and consequently,
\[
p_\ve^{-\alpha}=p^{-\alpha}\left (1-\alpha\ve \left ( \frac
{\Delta p}{p}+q\right)\right) \left(1 + O(\ve^{3/2})\right)
\]
Let
\[
k_\ve^{(\alpha)}(\x,\y)=\frac{k_\ve(\x,\y)}{p_\ve^{\alpha}(\x)p_\ve^{\alpha}(\y)}
\]
Then, the normalization factor $d_\ve^{(\alpha)}$ is given by
\[
d_\ve^{(\alpha)}(\x) = \int k_\ve^{(\alpha)}(\x,\y) p(\y)d\y =
p_\ve^{-\alpha}(\x) p^{1-\alpha}(\x)\left[1+\ve\left((1-\alpha)q
-\alpha\frac{\Delta p}{p} +\frac{\Delta
p^{1-\alpha}}{p^{1-\alpha}(\x)}\right)\right]
\]
Therefore, the asymptotic expansion of the backward operator gives
\[
T_b^{(\alpha)}\phi = \int_X
\frac{k^{(\alpha)}_\ve(\x,\y)}{d_\ve^{(\alpha)}(\x)}
\phi(\y)p(\y)d\y =\phi(\x) + \ve\left(\frac{\Delta(\phi
p^{1-\alpha})}{p^{1-\alpha}} - \phi \frac{\Delta
p^{1-\alpha}}{p^{1-\alpha}}\right)
\]
and its infinitesimal generator is
\[
\mathcal H_b \phi=\lim_{\ve\to 0}\frac{T_b-I}{\ve}\phi =
\frac{\Delta (\phi
p^{1-\alpha})}{p^{1-\alpha}}-\frac{\Delta(p^{1-\alpha})}{p^{1-\alpha}}\phi
\]
Inserting the expression $p=e^{-U}$ into the last equation gives
\[
\mathcal H_b \phi = \Delta \phi - 2 (1-\alpha) \nabla \phi \cdot
\nabla U
\]
Similarly, the form of the forward infinitesimal operator is
\[
\mathcal H_f \psi = \Delta \psi - 2 \alpha \nabla \psi \cdot
\nabla U + (2\alpha-1) \psi \left(\nabla U\cdot\nabla U - \Delta
U\right)
\]

\end{document}